\documentclass[11pt]{article}

\usepackage{amssymb,amsthm,amscd,array}

\let\ssection=\section
\renewcommand{\section}{\setcounter{equation}{0}\ssection}

\newcommand{\bbR}{\mathbb{R}}

\newcommand{\bbC}{\mathbb{C}}

\newcommand{\rD}{\mathrm{D}}

\newcommand{\ddiv}{\mathrm{div}}
\newcommand{\cD}{{\mathcal{D}}}

\newcommand{\Diff}{\mathrm{Diff}}

\newcommand{\rE}{\mathrm{E}}

\newcommand{\End}{\mathrm{End}}
\newcommand{\cF}{{\mathcal{F}}}

\newcommand{\gr}{\mathrm{gr}}

\newcommand{\rH}{\mathrm{H}}

\newcommand{\Hom}{\mathrm{Hom}}

\newcommand{\Id}{\mathrm{Id}}

\newcommand{\tr}{\mathrm{tr}}

\newcommand{\cS}{{\mathcal{S}}}

\newcommand{\Gl}{\mathrm{gl}}

\newcommand{\Sl}{\mathrm{sl}}

\newcommand{\Supp}{\mathrm{Supp}}

\newcommand{\Vect}{\mathrm{Vect}}

\newcommand{\half}{\frac{1}{2}}
\newcommand{\fg}{\mathfrak{g}}
\newcommand{\fh}{\mathfrak{h}}

\def\mod#1{\left|{#1}\right|}

\begin{document}



\def\d{\delta}
\def\g{\gamma}
\def\om{\omega}
\def\r{\rho}
\def\a{\alpha}
\def\b{\beta}
\def\s{\sigma}
\def\vfi{\varphi}
\def\l{\lambda}
\def\m{\mu}
\def\implies{\Rightarrow}

\oddsidemargin .1truein
\newtheorem{thm}{Theorem}[section]
\newtheorem{lem}[thm]{Lemma}
\newtheorem{cor}[thm]{Corollary}
\newtheorem{pro}[thm]{Proposition}
\newtheorem{ex}[thm]{Example}
\newtheorem{rmk}[thm]{Remark}
\newtheorem{defi}[thm]{Definition}


\title{Cohomology of the vector fields Lie algebra and modules of differential
operators on a smooth manifold}

\author{P.B.A. Lecomte
\thanks{Institute de Math\'ematiques,
 Universit\'e de Li\`ege, Sart Tilman,
Grande Traverse, 12 (B 37),
B-4000 Li\`ege,
BELGIUM, mailto:plecomte@ulg.ac.be
}
\and
 V.Yu. Ovsienko
\thanks{C.N.R.S., Centre de Physique Th\'eorique,
 Luminy -- Case 907,
F--13288 Marseille, Cedex 9,
FRANCE, mailto:ovsienko@cpt.univ-mrs.fr
}
}

\date{}

\maketitle

\thispagestyle{empty}

\begin{abstract}
Let $M$ be a smooth manifold, $\cS$ the space of polynomial on fibers
functions on $T^*M$ (i.e., of symmetric contravariant tensor fields).
We compute the first cohomology space of the Lie algebra, $\Vect(M)$, of 
vector fields on $M$ with coefficients in the space of linear differential operators
on $\cS$.  This cohomology space is closely related to the
$\Vect(M)$-modules, $\cD_\l(M)$, of linear differential operators on the space of
tensor densities on~$M$ of degree $\l$.
\end{abstract}


\section{Introduction and the Main Theorem}

Let $M$ be a smooth manifold and $\Vect(M)$ the Lie algebra of vector
fields on~$M$. 

The main purpose of this article is to study the cohomology of
$\Vect(M)$ with coefficients in the space of linear differential operators acting on
tensor fields. This cohomology is, actually, a natural generalization of the
Gelfand-Fuchs cohomology (i.e., of $\Vect(M)$-cohomology with coefficients in the
modules of tensor fields on $M$).

The problem of computation of such cohomology spaces naturally arises if one
considers \textit{deformations} of the $\Vect(M)$-module structure on
the space of tensor fields.

The general theory of deformations of Lie algebra modules is due to Nijenhuis
and Richardson \cite{NR,Ric}. Let $\fg$ be a Lie algebra and $V$ a
$\fg$-module, then the problem of deformation of the $\fg$-module structure on $V$
is related to the cohomology spaces: $\rH^1(\fg;\End(V))$ and
$\rH^2(\fg;\End(V))$. More precisely, the first cohomology
space classifies \textit{infinitesimal} deformation, while the second one
contains the obstructions to integrability of a given infinitesimal deformation.

The origin of our investigation is related to the space of scalar
linear differential operators on $M$ viewed as a module over $\Vect(M)$. It is
quite clear a-priori that this module should be considered as a deformation of the
corresponding module of \textit{symbols} (i.e., of polynomial on fibers
functions on $T^*M$). We are, therefore, led to study the first cohomology of
$\Vect(M)$ with coefficients in the
$\Vect(M)$-module of operators on the space of symbols.

\subsection{Differential operators on symmetric contravariant tensor
fields}\label{SkSpace}

Consider the space, $\cS(M)$ (or $\cS$ for short), of symmetric
contravariant tensor fields on $M$ (i.e., $\cS=\Gamma(STM)$). As
a $\Vect(M)$-module it is isomorphic to the space of smooth
functions on $T^*M$ polynomial on the fibers. Therefore, $\cS$ is a Poisson algebra
with a natural graduation given by the decomposition
\begin{equation}
\label{decomp}
\cS=\bigoplus_{k=0}^\infty\cS_k,
\end{equation}
where $\cS_k$ is the space of $k$-th order tensor fields. Obviously,
$\cS_0$ is isomorphic to $C^\infty(M)$ and $\cS_1$ to $\Vect(M)$. The Poisson bracket
on
$\cS$ is usually called the (symmetric) Schouten bracket (see e.g. \cite{Fuc}).

The action of $X\in\Vect(M)$ on $\cS$ is given by the Hamiltonian vector field
\begin{equation}
\label{Hamilton}
L_X=
\frac{\partial{}X}{\partial\xi_i}\,
\frac{\partial}{\partial{}x^i}
-
\frac{\partial{}X}{\partial{}x^i}\,
\frac{\partial}{\partial\xi_i}\,,
\end{equation}
where $(x,\xi)$ are local coordinates on $T^*M$ (we identified $X$ with the
first-order polynomial $X=X^i\xi_i$; the summation over repeated indices is
understood).

Let us introduce the space, $\cD(\cS)$, of all linear differential operators on
$\cS$. This space is a
$\Vect(M)$-module with a filtration
\begin{equation}
\label{filtr1}
\cD^0(\cS)
\subset
\cD^1(\cS)
\subset\cdots\subset
\cD^r(\cS)
\subset\cdots,
\end{equation}
where $\cD^r(\cS)$ is the space of $r$-th order differential
operators.

In this article we compute the first cohomology space
\begin{equation}
\label{DefCohom}
\rH^1(\Vect(M);\cD(\cS)).
\end{equation}
of $\Vect(M)$ acting on $\cD(\cS)$. 

Note that for $M=S^1$ this computation has been done in
\cite{LO1,BO} see also \cite{FF} for the case of the Lie algebra of formal vector
fields on $\bbR$.

\subsection{Modules of differential operators on tensor densities}

Let $\cF_\l(M)$ (or $\cF_\l$ in short) be the space of tensor densities of degree
$\l$ on~$M$ (i.e. the space of sections of the line bundle
$\Delta_\l(M)=\mod{\Lambda^nT^*M}^{\otimes\l}$ over $M$). Clearly,
$\cF_0\cong{}C^\infty(M)$ as a $\Vect(M)$-module, any two $\Vect(M)$-modules of
tensor densities are non-isomorphic (see also \cite{Fuc}).

Denote
$\cD_\l$ the space $\cD(\cF_\l)$ of linear differential operators on
$\cF_\l$. This space is an
associative (and, therefore, a Lie) algebra with the filtration by the order of
differentiation:
\begin{equation}
\label{filtr1}
\cD^0_\l\subset\cD^1_\l\subset\cdots\subset\cD^k_\l\subset\cdots
\end{equation}
The algebra $\cS$ is naturally identified with the associated graded algebra
$\gr(\cD_\l)$ that is, 
\begin{equation}
\label{principal}
\cD^k_\l/\cD^{k-1}_\l\cong\cS_k\,.
\end{equation}
The corresponding projection $\s_k:\cD^k_\l\to\cS_k$ is called the (principal)
\textit{symbol}.

\goodbreak

The associative algebra $\cD_\l$ can be naturally
interpreted as a non-trivial deformation of $\cS$ and constitutes one
of the main objects considered in \textit{deformation quantization}.

We will be interested, however, only in the $\Vect(M)$-module structure on $\cD_\l$
rather than in the whole associative (or Lie algebra) structure.
The (tautological) Lie algebra embedding $\Vect(M)\hookrightarrow\cD_\l$
\begin{equation}
\label{StandEmbed}
X\mapsto{}L^\l_X,
\end{equation}
where $L^\l_X$ is the Lie derivative on $\cF_\l$,
defines a $\Vect(M)$-module structure on $\cD_\l$.

\begin{rmk}
{\rm
If $M$ is oriented by a volume form $\Omega$,
then
\begin{equation}
\label{LieDerDens}
L^\l_X=L_X+\l\,\ddiv_\Omega{X}.
\end{equation}
Moreover, $\cD_\l$ and $\cD_\m$ are isomorphic associative algebras. However, as
$\Vect(M)$-modules they are isomorphic if and only if
$\l+\m=1$ \cite{DO,LMT}.
}
\end{rmk}

\subsection{The Main Theorem}

The space $\cD(\cS)$ is decomposed, as a
$\Vect(M)$-module, into the direct sum:
\begin{equation}
\label{Razlozhenie}
\cD(\cS)
=\bigoplus_{k,\ell}
\cD(\cS_k,\cS_\ell),
\end{equation}
where $\cD(\cS_k,\cS_\ell)\subset\Hom(\cS_k,\cS_\ell)$.
It would then suffice to compute the cohomology (\ref{DefCohom}) with coefficients
in each of these modules. Our main result is the following

\begin{thm}
\label{CohomThm}
If $\dim M\geq2$, then
\begin{equation}
\label{CohomResult}
\rH^1(\Vect(M);\cD(\cS_k,\cS_\ell))
=
\left\{
\begin{array}{lcl}
\bbR & , & \hbox{if} \quad k-\ell=2\\[6pt]
\bbR & , & \hbox{if} \quad k-\ell=1, \ell\neq0\\[6pt]
\bbR\oplus\rH^1_{\rm DR}(M) & , & \hbox{if} \quad k-\ell=0\\[6pt]
0 & , & \hbox{otherwise}
\end{array}
\right.
\end{equation}
where $\rH^1_{\rm DR}(M)$ is the first space of the de Rham cohomology of $M$.
\end{thm}

\noindent
The proof will be given in 
Section \ref{TheProof}.

\medskip

From now on we assume that $\dim{}M\geq2$.

\subsection{Differentiability}\label{peetre}

As a first step towards the proof of Theorem \ref{CohomThm},
we will prove now that any 1-cocycle on $\Vect(M)$ with values
in the space of differential operators $\cD(\cS_k,\cS_\ell)$ is locally
differentiable. Due to the well-known Peetre Theorem \cite{Pee}, this means that for
any $\g\in{\rm Z}^1(\Vect(M);\cD(\cS_k,\cS_\ell))$, the bilinear map
$(X,P)\mapsto\g(X)(P)$, where
$X\in\Vect(M)$ and $P\in\cS_k$, is local:
\begin{equation}
\label{Sup}
\Supp\,\g(X)(P)
\subset\Supp\,{X}\cap\,\Supp\,{P}
\end{equation}

\goodbreak

\begin{pro}
\label{PeetrePro}
Any 1-cocycle $\g$ on $\Vect(M)$ with values in $\cD(\cS_k,\cS_\ell)$ is local.
\end{pro}
\begin{proof}
Let $U\subset{}M$ be open and
$X\in\Vect(M)$ vanish on $U$. We have to show that
$
\g(X)_{|_U}=0.
$
Let $x_0$ be any point in $U$. As well-known, there exists a neighborhood
$V\subset{}U$ of $x_0$ and vector fields
$X_i,X_i'$, $i=1,\ldots,r$ on $V$ such that
$$
X=
\sum_{1\leq{}i\leq{}r}
[X_i,X_i']
$$
and 
$$
{X_i}
_{|_V}={X_i'}
_{|_V}=0,
$$
where $r$ depends only on the dimension of $M$. One has, using the fact that $\g$ is a
1-cocycle
$$
\g(X)_{|_V}=
\sum_{1\leq{}i\leq{}r}
\left(
L_{X_i}\g(X_i')_{|_V}-L_{X_i'}\g(X_i)_{|_V}
\right)=0.
$$
\end{proof}

\section{Non-trivial cohomology classes}\label{FormSection}

Let us now describe a natural basis of the above cohomology
spaces (\ref{CohomResult}).

\subsection{Case $k=\ell$}

Since $\Id\in\cD(\cS_k,\cS_k)$ is $\Vect(M)$-invariant, $c\mapsto{}c\,\Id$
maps any cocycle $c$ to a cocycle and thus induces a homomorphism
$\rH(\Vect(M);C^\infty(M))\to\rH(\Vect(M);\cD(\cS_k,\cS_k))$. Theorem~\ref{CohomThm}
states that it is an isomorphism in degree one.

Recall that $\rH(\Vect(M);C^\infty(M))$ is well-known (see \cite{Fuc}). In particular,
given a covariant derivation $\nabla$, the 1-cocycles are the maps
\begin{equation}
\label{OneCocycle}
c_{a,\omega}:
X
\mapsto
a\,\ddiv_\nabla(X)+i_X\omega,
\end{equation}
where $a\in\bbR$ and $\omega$ is a closed 1-form, $\ddiv_\nabla$ being the divergence
associated to $\nabla$. The cocycle (\ref{OneCocycle}) is a coboundary if and only if
$a=0$ and $\omega$ is exact.

\goodbreak

\subsection{Case $k=\ell+1$, $\ell\neq0$}\label{Konstrukciya1}

Consider the exact sequence of $\Vect(M)$-modules
\begin{equation}
\label{exact1}
\begin{CD}
0 @> >>\cD^{k-1}_\l @> >>\cD_\l^k @> >>\cS_k @> >>0.
\end{CD}
\end{equation}
Dividing out by $\cD_\l^{k-2}$ leads to the exact sequence
\begin{equation}
\label{exact1First}
\begin{CD}
0 @> >>\cS_{k-1} @> >>\cD_\l^k/\cD_\l^{k-2} @> >>\cS_k @> >>0
\end{CD}
\end{equation}
Assume $k\neq1$ and
$\l\neq1/2$. The sequence~(\ref{exact1First}) does not split \cite{LO}. Its cohomology
class is a non-zero element in $\rH^1(\Vect(M);\Hom(\cS_k,\cS_{k-1}))$ (see Appendix).
This class admits a representative with values in $\cD(\cS_k,\cS_{k-1})$,
since the $\Vect(M)$-actions in (\ref{exact1First}) are differential. It thus defines
a non-trivial class in $\rH^1(\Vect(M);\cD(\cS_k,\cS_{k-1}))$ which,
by Theorem~\ref{CohomThm}, is a basis of this space.



\subsection{Case $k=\ell+2$}\label{RelationSec}

If $\l=1/2$, it is shown in \cite{LO} that the sequence~(\ref{exact1First}) is split
and that the sequence
\begin{equation}
\label{exact2}
\begin{CD}
0@> >>\cS_{k-2}@> >>
\cD_{1/2}^{k}/\cD_{1/2}^{k-3}@> >>
\cD_{1/2}^{k}/\cD_{1/2}^{k-2}@> >>0
\end{CD}
\end{equation}
is not. Moreover, the splitting of (\ref{exact1First}) is given by differential
projectors.
Since (\ref{exact1First}) is split, the class
$[\cD^{k-1}_{1/2},\cD_{1/2}^k]$ of (\ref{exact1}) belongs to
$\rH^1(\Vect(M);\Hom(\cS_{k},\cD_{1/2}^{k-2}))$. Since (\ref{exact2}) is not split, its
projection
${\s_{k-2}}_\sharp\,[\cD^{k-1}_{1/2},\cD_{1/2}^k]$ is non-zero (see Lemma
\ref{SequenceLemma} from Appendix). 

As in the previous case, this projection is easily seen to admit a representative with
values in $\cD(\cS_k,\cS_{k-2})$. Hence, it provides a basis of
$\rH^1(\Vect(M);\cD(\cS_k,\cS_{k-2}))$.

\begin{rmk}
\label{RemRes}
{\rm
In the above Subsections \ref{Konstrukciya1} and \ref{RelationSec} we have associated
non-trivial cohomology classes to the exact sequences (\ref{exact1First}) and
(\ref{exact2}). It is important to note that these classes are ``natural'' in the
following sense. For any open subset $U\subset{}M$ their restrictions to $U$ are
precisely the classes associated to the same sequences upon $U$. }
\end{rmk}

\section{Projectively equivariant cohomology}\label{ProjectCohom}

Throughout this section we put $M\cong\bbR^n$ and $n\geq2$.

\subsection{The Lie algebra of infinitesimal projective
transformations}\label{Introducing}

The main idea of our proof of Theorem~\ref{CohomThm} is to use the
filtration with respect to the Lie subalgebra
\begin{equation} 
\label{Definitionsl2}
\Sl(n+1,\bbR)\subset\Vect(\bbR^n).
\end{equation} 
It is suggested by the fact that the exact sequence (\ref{exact1}) that generate our
cohomology is split as a sequence of $\Sl(n+1,\bbR)$-modules \cite{LO}. 
In some sense, this Lie subalgebra plays the
same r\^ole in our approach as the linear subalgebra $\Gl(n,\bbR)$ in the
traditional one (cf.~\cite{Fuc}).

Recall that the standard action of the Lie
algebra $\Sl(n+1,\bbR)$ on $\bbR^n$ is generated by the vector fields
\begin{equation} 
\label{sl2}
X_i=
\frac{\partial}{\partial{}x^i}\,, 
\qquad 
X_{ij}=
x^i\frac{\partial}{\partial{}x^j}\,, 
\qquad
\bar{X}_i=
x^i{\cal E}\,,
\end{equation}
where
\begin{equation} 
{\cal E} = x^i\frac{\partial}{\partial{}x^i}\,.
\label{Euler}
\end{equation}
Observe in particular that $X_i$ and $X_{ij}$ generate an action of the Lie
algebra $\Gl(n,\bbR)\ltimes\bbR^n$.

\subsection{Computing the relative cohomology space}

In this section we will compute the first space of the so-called relative
cohomology of $\Vect(\bbR^n)$,
i.e. the cohomology of the complex of $\Vect(\bbR^n)$-cochains vanishing on the
subalgebra $\Sl(n+1,\bbR)$. We will prove the following

\goodbreak

\begin{thm}
\label{RelativeTheorem}
If $n\geq2$, then
\begin{equation}
\label{RelatCohom}
\rH^1(\Vect(\bbR^n),\Sl(n+1,\bbR);\cD(S_k,\cS_\ell))=
\left\{
\begin{array}{lcl}
\bbR & , & \hbox{if} \quad k-\ell=2\\[6pt]
\bbR & , & \hbox{if} \quad k-\ell=1, \ell\neq0\\[6pt]
0 & , & \hbox{otherwise}
\end{array}
\right.
\end{equation}
\end{thm}

\subsection{Equivariance property}\label{EquivSection}

We begin the proof with a simple observation.

Let $\fh\subset\fg$ be a Lie subalgebra and $V$ a $\fg$-module.
If $c:\fg\to{}V$ is a 1-cocycle
such that $c_{|_\fh}\equiv0$, then it is \textit{equivariant} with respect to $\fh$
i.e.
\begin{equation}
\label{equivCoc}
L_X(c(Y))=
c([X,Y]),
\qquad X\in\fh,
\end{equation}
where $L$ stays for the $\fg$-action on the module $V$.

Consequently, our strategy to compute the space of relative cohomology
(\ref{RelatCohom}) consists, first, in classifying
the $\Sl(n+1,\bbR)$-equivariant linear maps $c:\Vect(\bbR)\to\cD(\cS_k,\cS_\ell)$
vanishing on $\Sl(n+1,\bbR)$ and, second, to isolate among them the 1-cocycles.

\subsection{Commutant of the affine Lie algebra}

Consider the space of polynomials
$\bbC[x,\xi]=\bbC[x^1,\ldots,x^n,\xi_1,\ldots,\xi_n]$ as a submodule of $\cS$ under
the action of $\Sl(n+1,\bbR)$. We need to compute the \textit{commutant}
of the subalgebra $\Gl(n,\bbR)\ltimes\bbR^n$, i.e. the algebra of
differential operators on
$\bbC[x,\xi]$ commuting with the $\Gl(n,\bbR)\ltimes\bbR^n$-action.

The differential operators on $\bbC[x,\xi]$ given by
\begin{equation}
\label{Div}
\rE=\xi_i\,
\frac{\partial}{\partial\xi_i},\qquad
\rD=\frac{\partial}{\partial x^i}\,
\frac{\partial}{\partial\xi_i}
\end{equation}
commute with the $\Gl(n,\bbR)\ltimes\bbR^n$-action.
Let us recall the classical result of the Weyl
invariant theory (see \cite{Wey}).
\begin{pro}
\label{WeylPro}
The algebra of differential operators on
$\bbC[x,\xi]$ commuting with the action of the
affine Lie algebra, is generated by $\rE$ and $\rD$.
\end{pro}
\noindent
We will call the operators (\ref{Div}) the Euler
operator and the divergence operator respectively.
The eigenspaces of $\rE$ are obviously consist of homogeneous
polynomials in $\xi$.

\begin{cor}
\label{NiceCor}
The operator $\rD^{k-\ell}$ is the unique (up to a constant)
$\Gl(n,\bbR)\ltimes\bbR^n$-equiva\-riant differential operator from $\cS_k$ to
$\cS_\ell$.
\end{cor}
\begin{proof}
Any differential operator on $\cS_k$ is indeed determined by its values on the
subspace~$\bbC[x,\xi]$.
\end{proof}

The Euler operator $\rE$ is clearly equivariant with respect to the whole
$\Vect(\bbR^n)$. We will need the commutation relations of the operator $\rD$
with the quadratic generators of~$\Sl(n+1,\bbR)$.
\begin{lem}
For $\bar{X}_i$ as in (\ref{sl2}), one has
\begin{equation}
[L_{\bar{X}_i},\rD]
=
\Big(2\rE+(n+1)\Big)\circ\frac{\partial}{\partial{\xi_i}}\,,
\label{relation}
\end{equation}
\end{lem}
\begin{proof}:
straightforward.
\end{proof}

\subsection{Bilinear $\Gl(n,\bbR)\ltimes\bbR^n$-invariant operators}

We also need to classify the bilinear $\Gl(n,\bbR)\ltimes\bbR^n$-invariant differential
operators. For that purpose, let us use a natural identification
\begin{equation}
\label{IdentBin}
\bbC[x,\xi]\otimes\bbC[y,\eta]
\cong
\bbC[x,\xi,y,\eta].
\end{equation}
There are, obviously, four invariant differential operators
$\rD_{(x,\xi)},\rD_{(y,\eta)}$ (the divergence
operators with respect to the first and the second
arguments)
and $\rD_{(x,\eta)},\rD_{(y,\xi)}$ (the operators of contraction in terms of
tensors). Applying again \cite{Wey} one gets the following
\begin{pro}
\label{AffPro}
Every bilinear differential operator
\begin{equation}
\label{Bilin1}
\cS_j\otimes\cS_k\to\cS_\ell
\end{equation}
invariant with respect to the action of the affine Lie algebra, is a homogeneous
polynomial in $\rD_{(x,\xi)},\rD_{(x,\eta)},\rD_{(y,\xi)}$ and $\rD_{(y,\eta)}$
of degree $j+k-\ell$.
\end{pro}

We are now ready to start the proof of Theorem \ref{RelativeTheorem}.

\subsection{Bilinear $\Sl(n+1,\bbR)$-equivariant operators}

In view of Section \ref{EquivSection}, we will now classify
the $\Sl(n+1,\bbR)$-equivariant linear differential maps
\begin{equation}
\label{LinearEquiv}
c:\Vect(\bbR^n)\to\cD(\cS_k,\cS_\ell)
\end{equation}
vanishing on the subalgebra $\Sl(n+1,\bbR)\subset\Vect(\bbR^n)$.
We can, equivalently, consider the equivariant bilinear maps
\begin{equation}
\label{Bilin}
C:\cS_1\otimes\cS_k\to\cS_{k-p},
\end{equation}
where $p=k-\ell$.

By Proposition \ref{AffPro}, any such operator is of the
form
\begin{equation}
\label{Ansatz}
\begin{array}{rcl}
C&=&
\displaystyle\sum_{s=0}^{p+1}
\left(
\frac{\a_s}{s!(p-s+1)!}\,
\rD_{(x,\eta)}^{s}\rD_{(y,\eta)}^{p-s+1}\right.\\[12pt]
&&\qquad
\left.\displaystyle+\frac{\b_s}{(s-1)!(p-s+1)!}\,
\rD_{(x,\xi)}\rD_{(x,\eta)}^{s-1}\rD_{(y,\eta)}^{p-s+1}\right)_{
\Big|\!\begin{array}{l}
y=x \\
\eta= \xi
\end{array}
}\\[12pt] &&
{\displaystyle+\sum_{s=0}^p
\frac{\g_s}{s!(p-s)!}\,
\rD_{(y,\xi)}\rD_{(x,\eta)}^s\rD_{(y,\eta)}^{p-s}}_{\Big|\!
\begin{array}{l}
y=x \\
\eta= \xi
\end{array}}
\end{array}
\end{equation}
where $\a_s,\b_s,\g_s\in\bbR$. 

Moreover,
\begin{equation}
\label{NulForSmall}
\a_s=\b_s=\g_s=0
\qquad
\hbox{for}
\qquad
s<2,
\end{equation}
since $C$ vanishes on the affine subalgebra and
\begin{equation}
\label{Zero}
(k-p)\a_2+(n+1)\b_2+(p-1)\g_2=0,
\end{equation}
since $C$ vanishes on the quadratic generators $\bar{X}_i$ of $\Sl(n+1,\bbR)$.

For $k=p$, we have not to take into account the coefficients $\a_s$ in the expression
(\ref{Ansatz}) because the corresponding terms vanish when applied to
$\cS_1\otimes\cS_k$.

It is quite easy, using (\ref{relation}) and
analogous relations with the operators $\rD_{(x,\eta)},\rD_{(y,\xi)}$ and
$\rD_{(y,\eta)}$, to obtain the necessary and sufficient condition for the coefficients
in (\ref{Ansatz}) for~$C$ to be equivariant. One gets the following recurrence
relations:
\begin{eqnarray}
(s-1)\,\a_{s+1}-(2k+n-p+s-1)\,\a_s-\g_s
&=&0
\label{sys1}
\\
(s-1)\,\b_{s+1}-(2k+n-p+s-1)\,\b_s-\g_s
&=&0
\label{sys2}
\\
(s-2)\,\g_s-(2k+n-p+s-1)\,\g_{s-1}
&=&0
\label{sys3}
\\
(k-p)\,\a_{s+1}+(n+1)\,\b_{s+1}+(p-s)\,\g_{s+1}+(k-p+s)\,\g_s
&=&0
\label{sys4}
\end{eqnarray}
where $2\leq{}s\leq{}p$.
(For $k=p$, equation (\ref{sys1}) has not to be taken into account.)

Now, to solve the system (\ref{Zero}-\ref{sys4}), we need the following technical
\begin{lem}
\label{Simplify}
If $\a_s,\b_s,\g_c$ verify the equations (\ref{sys1}-\ref{sys3}) and (\ref{Zero}),
then $\a_s,\b_s,\g_c$ verify the equation (\ref{sys4}).
\end{lem}
(A similar result holds true when $k=p$.)
\begin{proof}
Check that for $s=1$, the equation (\ref{sys4}) coincides with (\ref{Zero}),
the result follows then by induction.
\end{proof}

It is now very easy to get the complete solution of the system
(\ref{Zero}-\ref{sys3}). One has the following four cases.

(a) For $p=0$ and for $(p=1,k=1)$ there is no solution.

\medskip

(b) For $(p=1,k\geq2)$ the system has a one-dimensional space of
solutions spanned by
\begin{equation}
\label{FirstCocycle}
C_1=\half\,\rD^2_{(x,\eta)}+
\frac{k-1}{n+1}\,\rD_{(x,\xi)}\rD_{(x,\eta)}\,,
\end{equation}
which is, in fact, a just a solution of the equation (\ref{Zero}).
\begin{rmk}
{\rm
One readily checks that the operator $c(X)\in\cD(\cS_k,\cS_{k-1})$ given by
(\ref{FirstCocycle}) coincides (up to a constant) with the operator of contraction
with the tensor field
\begin{equation}
\label{FirstCocycleTensor}
c_1(X)=
\left(
\frac{\partial}{\partial{}x^i}
\frac{\partial}{\partial{}x^j}(X^\ell)
+
\frac{2}{n+1}\,
\d_j^\ell
\frac{\partial}{\partial{}x^i}
\frac{\partial}{\partial{}x^s}(X^s)
\right)
dx^i{}dx^j
\otimes\xi_\ell
\end{equation}
This expression is obviously a 1-cocycle. The expression (\ref{FirstCocycleTensor})
is known in the literature as the Lie derivative of a flat 
\textit{projective connection} (cf., e.g., \cite{kh}).
}
\end{rmk}

(c) For $p\geq2,k>p$, the system (\ref{sys1}-\ref{sys3}) under the
condition (\ref{Zero}), has a two-dimensional space of solutions parametrized by
$(\a_2,\b_2)$.

\medskip

(d) For $p=k\geq2$, the equation (\ref{sys1}) should be discarded.
The space of solutions is again
one-dimensional.

\subsection{Projectively invariant cocycles}

We will now determine which of the $\Sl(n+1,\bbR)$-equivariant maps
(\ref{LinearEquiv}) classified in the preceding section are 1-cocycles. Let us examine
separately the cases~(b)-(d).

\medskip

(b) In the simplest case, $p=1$, one easily checks that the unique
$\Sl(n+1,\bbR)$-equivariant map (\ref{FirstCocycle}),
indeed, defines a 1-cocycle on $\Vect(\bbR^n)$ with values
in $\cD^0(\cS_k,\cS_{k-1})$.

\medskip

(c) 
The cocycle relation adds the equation $\b_3=2\b_2$ to the general
system (\ref{sys1}-\ref{sys3}).

In the case $p=2$, one checks by a straightforward computation, that the solutions
are the constant multiples of the solution given by
\begin{equation}
\label{CocSecondDef}
\begin{array}{rcl}
\a_2 &=& 2,\\[4pt]
\a_3 &=& 2k+n+1,\\[4pt]
\b_2 &=& 1, \\[4pt]
\b_3 &=& 2, \\[4pt]
\d_2 &=& \!\! -(2k+n-3).
\end{array}
\end{equation}

In the case $p>2$, the only solution of the system (\ref{sys1}-\ref{sys3}) together
with the equation $\b_3=2\b_2$ is zero.

\medskip

(d) If $k=p$, then the non-trivial solutions of the system are cocycles if and only if
$k=p=2$. This cocycle is precisely of the form (\ref{CocSecondDef}) disregarding
$\a_2$ and
$\a_3$.

\begin{pro}
\label{NontrivPro}
The 1-cocycles on $\Vect(\bbR^n)$ defined by the formul{\ae} (\ref{FirstCocycle}) and
(\ref{CocSecondDef}) are non-trivial.
\end{pro}
\begin{proof}
This follows immediately from Sections \ref{Konstrukciya1}, \ref{RelationSec} and the
fact that the sequence (\ref{exact1}) is split when restricted to $\Sl(n+1,\bbR)$, see
\cite{LO}. Let us also give an elementary proof.

Recall that a 1-cocycle on $\Vect(\bbR^n)$ with values in
$\cD(\cS_k,\cS_\ell)$ is a coboundary if it is of the form
$X\mapsto{}[L_X,B]$ for some
$B\in\cD(\cS_k,\cS_\ell)$. Moreover, the 1-cocycle vanishes on
$\Sl(n+1,\bbR)$ if and only if $B$ is $\Sl(n+1,\bbR)$-equivariant.
\begin{lem} (cf. \cite{Lec1}).
\label{LemInv}
If $k\neq\ell$, there is no $\Sl(n+1,\bbR)$-equivariant operators
$B\in\cD(\cS_k,\cS_\ell)$ different from zero.
\end{lem}
\begin{proof}
In virtue of Corollary~\ref{NiceCor}, the property of $\Sl(n+1,\bbR)$-equivariance
implies, in particular, that $B$ has to be proportional to $\rD^{k-\ell}$. Now,
the commutation relation (\ref{relation}) shows that this operator can never be
$\Sl(n+1,\bbR)$-equivariant.
\end{proof}
Proposition \ref{NontrivPro} follows.
\end{proof}

\subsection{Proof of Theorem \ref{RelativeTheorem}}

We have shown that there
exist unique (up to a constant) 1-cocycles $c_1$ and $c_2$ on $\Vect(\bbR^n)$ with
values in $\cD(\cS_k,\cS_{k-1})$ and
$\cD(\cS_k,\cS_{k-2})$ respectively, vanishing on $\Sl(n+1,\bbR)$. These cocycles
define non-trivial classes of relative cohomology.

\medskip

Theorem \ref{RelativeTheorem} is proven.

\section{Proof of Theorem \ref{CohomThm}}\label{TheProof}

Using the filtration with respect to the
subalgebra $\Sl(n+1,\bbR)$, we will first prove Theorem~\ref{CohomThm} in the case
when $M$ is a vector space and then extend it to an arbitrary manifold.
To that end, we need some more information about the cohomology of $\Sl(n+1,\bbR)$.

\subsection{Cohomology of $\Sl(n+1,\bbR)$}

The cohomology of the Lie algebra $\Sl(n+1,\bbR)$ with coefficients in
$\cD(\cS_k,\cS_\ell)$ has been computed in \cite{Lec1}.
\begin{thm}
\label{SlThm}
The space of cohomology $\rH(\Sl(n+1,\bbR);\cD(\cS_k,\cS_\ell))$ is trivial for
$k\neq\ell$, for $k=\ell$ it is isomorphic to the Grassman algebra of invariant
functionals on
$\Gl(n,\bbR)$:
\begin{equation}
\label{CohSl}
\rH(\Sl(n+1,\bbR);\cD(\cS_k,\cS_k))
=
\left(
\bigwedge\textstyle\Gl(n,\bbR)^*
\right)^{\Gl(n,\bbR)}
\end{equation}
\end{thm}
In particular, 
\begin{equation}
\label{CohSlParticul}
\rH^1(\Sl(n+1,\bbR);\cD(\cS_k,\cS_\ell))=
\left\{
\begin{array}{rl}
\bbR,& k=\ell \\
0,& \hbox{otherwise}
\end{array}
\right.
\end{equation}
and the class of the 1-cocycle $X\mapsto\ddiv(X)\Id$ spans that space in the
case $k=\ell$.
(In fact, it
corresponds to the invariant function $\tr:\Gl(n,\bbR)\to\bbR$.) 
Note that this cocycle is just the restriction to $\Sl(n+1,\bbR)$ of the cocycle
$c_{1,0}$, see (\ref{OneCocycle}).

\subsection{The case of $\bbR^n$}\label{VectorSection}

The restriction of a 1-cocycle $c:\Vect(\bbR^n)\to\cD(\cS_k,\cS_\ell)$ to
$\Sl(n+1,\bbR)$ is a 1-cocycle on $\Sl(n+1,\bbR)$. If $k\neq\ell$, then this
restriction is trivial and, therefore, $c$ is cohomological to a 1-cocycle on
$\Vect(\bbR^n)$ vanishing on $\Sl(n+1,\bbR)$; if $k=\ell$, then the
restriction of $c$ to $\Sl(n+1,\bbR)$ is cohomological to
$c_{1,0}$ and so $c-c_{1,0}$ is, again, cohomological to a 1-cocycle on
$\Vect(\bbR^n)$ vanishing on $\Sl(n+1,\bbR)$.
The result then follows from Theorem~\ref{RelativeTheorem}.

Theorem \ref{CohomThm} is proven for the special case $M=\bbR^n$.

\subsection{The general case}

Let us now prove Theorem \ref{CohomThm} for an arbitrary manifold $M$. 
Consider a 1-cocycle $c$ on $\Vect(M)$ with values in $\cD(\cS_k,\cS_\ell)$.

\medskip

(a) If $k-\ell\neq0,1,2$, then in any domain of chart $U\cong\bbR^n$, the restriction
$c_{|_U}$ is a coboundary, that is $c(X)_{|_U}=L_X(S_U),$ where
$S_U\in\cD(\cS_k,\cS_\ell)$ is some operator on $U$. But, on $U\cap{}V$, one has
$c(X)_{|_{U\cap{}V}}=L_X(S_U)=L_X(S_V)$ and so the operator $S_U-S_V$ is invariant.
Lemma \ref{LemInv} implies $S_U-S_V=0$. Therefore, the $S_U$'s are the restrictions of
some globally defined $S\in\cD(\cS_k,\cS_\ell)$ and $c$ is its coboundary.

\medskip

(b) If $k-\ell=1$ or $2$, 
it follows from Theorem \ref{CohomThm} for $M=\bbR^n$ that 
the class of $c_{|_U}$ is determined up to a constant. In view of Remark \ref{RemRes},
one has thus
\begin{equation}
\label{LocalEquation}
c_{|_U}
=
\a_U\,\g_{|_U}+L_X(S_U),
\end{equation}
for some $\a_U\in \bbR$ and $S_U$ as above, where $\g$ is a representative of one of
the classes associated to the sequences (\ref{exact1First}) and (\ref{exact2})
respectively.
On $U\cap{}V$ one obviously has $\a_U=\a_V$ and $S_U=S_V$ since 
\[
\left(
\a_U-\a_V
\right)
\g_{|_{U\cap{}V}}
=
\partial
\left(
S_U-S_V
\right),
\]
$\g_{|_{U\cap{}V}}$ is non-trivial and, as above, $S_U-S_V$
is invariant.

\medskip

(c) If $k-\ell=0$, 
one has
\begin{equation}
\label{LocalEquationBis}
c_{|_U}
=
\a_U\,{c_{1,0}}_{|_U}+L_X(S_U).
\end{equation}
Once again, $\a_U=\a_V\,(:=a)$ and $S_U-S_V$
is invariant, but 
any invariant operator in $\cD(\cS_k,\cS_k)$
is proportional to the identity so that
$S_U-S_V=\b_{UV}\,\Id$, where $\b_{UV}$ is a constant. It is clear that the
$\b_{UV}$'s define a $\check{\rm C}$ech 1-cocycle. If now $\omega$ is a closed 1-form
representing the corresponding de Rham class, one easily sees that
$c$ is cohomologous to~$c_{a,\omega}$.

\medskip

Theorem \ref{CohomThm} is proven.

\section{Cocycles associated to a connection}\label{ConnectCocycles}

Using a torsion free covariant derivation $\nabla$, it is possible to construct
globally defined cocycles spanning $\rH^1(\Vect(M);\cD(\cS_k,\cS_\ell))$ for
$k-\ell=1,2$.

\subsection{Lie derivative of a connection}\label{Konstrukciya1}

For each vector field $X$, the Lie derivative
\[
L_X(\nabla) : (Y,Z) \mapsto [X,\nabla_Y Z] -
\nabla_{[X,Y]} Z-\nabla_Y[X,Z]
\]
of $\nabla$ is well-known to be a
symmetric $(1,2)$-tensor field. It yields a non-trivial
1-cocycle
\[
X
\mapsto
L_X(\nabla)
\]
on~$\Vect(M)$ with values in $\Gamma(\bigotimes_2^1TM)$.
Therefore, for $k\geq2$, the contraction
\begin{equation}
\label{TheFirstCoc1}
\g^\nabla_1(X)(P)
=
\langle
P,L_X(\nabla)
\rangle,
\qquad
P\in\cS_k,
\end{equation}
defines a 1-cocycle on~$\Vect(M)$ with values in $\cD^0(\cS_k,\cS_{k-1})$.

\subsection{Second-order cohomology class and the Vey cocycle}\label{Konstrukciya2}

The last case, $\ell=k-2$, is
directly related to deformation quantization.

For any symplectic manifold $V$, there exists a non-zero class
in $\rH^2(C^\infty(V);C^\infty(V))$. It is given by so-called \textit{Vey cocycle}
usually denoted $S^3_\Gamma$ (see \cite{WL} and \cite{Rog} for explicit
construction using a connection $\Gamma$ on $V$).

In the particular, if $V=T^*M$ one can choose the connection so that $S^3_\Gamma$ is
homogeneous of weight $-3$, namely, restricted to $\cS\subset{}C^\infty(T^*M)$,
\begin{equation}
\label{S3}
S^3_\Gamma:
\cS_k\otimes\cS_\ell\to\cS_{k+\ell-3},
\end{equation}
see \cite{WL1} (e.g. choosing $\Gamma$ as a lift of $\nabla$ to $T^*M$).
It follows easily from (\ref{S3}) that the map $\Vect(M)\to\cD(\cS_k,\cS_{k-2})$
defined by
\begin{equation}
\label{VeyCoc1}
\g_2^\nabla(X)(P)=
S^3_\Gamma(X,P),
\qquad
P\in\cS_k.
\end{equation}
is a 1-cocycle.

\section{Appendix: approximations of the class of a short exact
sequence of modules}\label{Appendix1}

\subsection{Class of a short exact sequence of $\fg$-modules}\label{ShortExact}

We will need some general information about short exact sequences of filtered
modules.

Let $\fg$ be a Lie algebra. Consider an exact sequence of $\fg$-modules
\begin{equation}
\label{s1}
\begin{CD}
0 @> >>A @>{i}>> B @>{j}>> C @> >>0
\end{CD}
\end{equation}
It is characterized by an element of
$\rH^1(\fg;\Hom(C,A))$ (cf.~\cite{Fuc}, Sec. 1.4.5). It will be convenient to denote
it $[A,B]$.
Recall that
if $\tau : C \to B$ is a section of $j$, then $[A,B]$ is the class of the
1-cocycle
$\g^\tau: \fg\to\Hom(C,A)$  given by
\begin{equation}
\label{AbstractCoc}
\g^\tau(X)(T) 
= 
i^{-1}(X.\tau(T) - \tau(X.T)),
\end{equation}
where $X\in\fg$ and $T\in{}C$ (this expression is well defined since $X.\tau(T) -
\tau(X.T)\in\ker{j}$).

Given a submodule $V$ of $A$, one has the following commutative diagram:
\begin{equation}
\label{SubmodDiagramm}
\begin{CD}
@. 0 \\
@. @VV V \\
@. V \\
@.@VV{i_V}V\\
0 @> >>A @> >> B @> >> C @> >>0 \\
@.@VV{\pi_A}V @VV{\pi_B}V @VV{{\rm Id}}V \\
0 @> >>A/V @> >> B/V @> >> C @> >>0 \\
@. @VV V \\
@. 0
\end{CD}
\end{equation}
where $i_V$ is the injection of $V$ into $A$ and $\pi_A$, $\pi_B$ are the
projections. 

\goodbreak

One has the relation $[A/V,B/V] = {\pi_A}_\sharp\,[A,B]$.
Moreover, the left vertical of 
(\ref{SubmodDiagramm})
leads  to the exact
triangle
\begin{equation}
\label{TriangleDiagramm}
\begin{array}{cl}
\mathrm{H}(\fg;\mathrm{Hom}(C,V)) \\
& \searrow {i_V}_\sharp \\
\Delta\;\big\uparrow & \mathrm{H}(\fg;\mathrm{Hom}(C,A)) \\
& \swarrow {\pi_A}_\sharp \\
\mathrm{H}(\fg;\mathrm{Hom}(C,A/V))
\end{array}
\end{equation}
where $\Delta$ is the connecting homomorphism. One easily obtains the following
\begin{pro}
\label{SequencePro}
\hfill\break
(i) The class $[A/V,B/V]$ vanishes if and only if 
$[A,B]\in\mathrm{im}\,{i_V}_\sharp$. 
\hfill\break
(ii) If $[A/V,B/V]=0$ then the class of the
exact sequence
\begin{equation}
\label{ExactB}
\begin{CD}
0 @> >>V @>{i \circ i_V} >> B @> >> B/V @> >>0
\end{CD}
\end{equation}
is $[V,B] = [V,A] + [A,B]$
and vanishes if and only if $[V,A]=[A,B]=0$.
\end{pro}

\subsection{Case of a filtered module}\label{Filtr}

Consider now a flag of filtered $\fg$-modules 
$A_0\subset{}A_1\subset\cdots\subset{}A_r\subset\cdots$ and put $S_r =A_r/A_{r-1}$.
Let us study the classes $[A_r,A_{r+1}]$ of the sequences
\begin{equation}
\label{BolshayaSeq}
\begin{CD}
0 @> >>A_r @> >> A_{r+1} @> >> S_{r+1} @> >>0.
\end{CD}
\end{equation}
The quotient by $V= A_{r-1}$, leads to its ``first approximation'':
\begin{equation}
\label{s3}
\begin{CD}
0 @> >>S_{r} @> >> A_{r+1}/A_{r-1} @> >> S_{r+1} @> >>0
\end{CD}
\end{equation}
If the sequence (\ref{s3})
is split, then $[A_r,A_{r+1}]\in\rH^1(\fg;\Hom(S_{r+1},A_{r-1}))$
and so
\begin{equation}
\label{SplitClass}
[A_{r-1},A_{r+1}] = [A_{r-1},A_r]+[A_r,A_{r+1}]
\end{equation}
by Proposition \ref{SequencePro}.

The next approximation is a result of the quotient by $A_{r-2}$.
Let $\pi_r:A_r\to{}A_r/A_{r-1}$ be the projection to the
quotient-module.
\begin{lem}
\label{SequenceLemma}
If the sequence (\ref{s3}) is split for all $r>0$, but the sequences
\begin{equation}
\label{SecondApprox}
\begin{CD}
0 @> >>S_{r-1} @> >> A_{r+1}/A_{r-2} @> >>A_{r+1}/A_{r-1} @> >>0,
\end{CD}
\end{equation}
for $r>1$ are not split, then the class
${\pi_{r-1}}_\sharp[A_r,A_{r+1}]$
does not vanish.
\end{lem}
\begin{proof}
Since the sequence (\ref{s3}) is split, one has
${\pi_{r-1}}_\sharp[A_{r-1},A_r]=0$. 
If in addition ${\pi_{r-1}}_\sharp[A_r,A_{r+1}]
= 0$, then by Proposition  \ref{SequencePro}
\[
[A_{r-1}/A_{r-2},A_{r+1}/A_{r-2}]=
{\pi_{r-1}}_\sharp
[A_{r-1},A_{r+1}]=
{\pi_{r-1}}_\sharp
\left(
[A_{r-1},A_r]
+[A_r,A_{r+1}]
\right)
= 0
\]
and the sequence (\ref{SecondApprox}) is split.
\end{proof}

\medskip
\noindent
{\it Acknowledgments}. It is a pleasure to acknowledge
numerous fruitful discussions
with  C.~Duval. We are also thankful to M.~De~Wilde, V.~Fock and C.~Roger for helpful
suggestions.

\newpage

\end{document}